\documentclass[12pt,hidelinks]{article}

\title{\(\lf\): a Foundational Higher-Order Logic}
\author{Zachary Goodsell and Juhani Yli-Vakkuri}

\input{Preamble.sty}

\begin{document}

\maketitle

\section{Introduction}

This paper presents\footnote{The system has already been discussed in the literature: see \cite{Williamson2023-BOP}: 220, which cites unpublished work by the authors. However, no explicit formulation of the system, \textit{qua} formal system, has been published to date.} a new system of logic, \(\lf\), that is intended to be used as the foundation of the formalization of science. That is, deductive validity according to \(\lf\)\footnote{And thus according to any useful conservative extensions thereof, of which we highlight $\lf_\iota$ and $\lf_\varepsilon$ in this paper (\S\ref{Sec-IotaAndEpsilon}).} is to be used as the criterion for assessing what follows from the verdicts, hypotheses, or conjectures of any science. In work currently in progress, we argue for the unique suitability of \(\lf\) for the formalization of logic, mathematics, syntax, and semantics. The present document specifies the language and rules of \(\lf\), lays out some key notational conventions, and states some basic technical facts about the system.

\(\lf\) is a system of higher-order logic based on that of \textcite{Church1940-CHUAFO}, and is closely related to the now-standard system of \textcite{Henkin1950-HENCIT}. \(\lf\) improves on these systems by being \textit{intensional}, like Henkin's system, but not \textit{extensional}, unlike Henkin's. An intensional system is one in which provably equivalent formulae may be substituted in any context. Intensionality suffices for the provability of all the equalities we need to be able to prove in mathematics, which we therefore also need to be able to prove throughout the rest of science, such as \(1=0+1\). In Church's system, \(1=0+1\) cannot be proved, provided that the numerals are defined, following \textit{Principia Mathematica} (\cite{Whitehead1910-WHIPM-8}), so as to be suitable for counting.

Henkin's system is an intensional extension of Church's system and is in that respect an improvement on it, but it incorporates the problematic assumption of \textit{extensionality}, which says that the material equivalence of propositions (type $t$\footnote{Or type $o$ in Church's and Henkin's terms.}) is sufficient for their identity. Extensionality is problematic for a variety of reasons, most obviously in applications of the theory of probability. Suppose that a fair coin is to be tossed, and let \(H\) and \(T\) be sentences that respectively formalize `the coin lands heads' and `the coin lands tails'. It would be normal to assume each of
\begin{align*}
    \mathrm{Pr}(H\wedge T)=0\quad \quad \mathrm{Pr}(H)=0.5\quad\quad \mathrm{Pr}(H\vee T)=1,
\end{align*}
but this trio is simply inconsistent in an extensional system. The derivation of a contradiction from the trio in an extensional system requires no extralogical assumptions about probability function \(\mathrm{Pr}\) or about anything else;\footnote{We assume that, like the numerals, the real number terms are treated as abbreviations of some suitable purely logical terms, so that \(0\neq 0.5\), \(0\neq 1\), and \(0.5\neq 1\) all abbreviate theorems.} the point is that in an extensional system we can prove that there are only two propositions, so no function can give them more than two values.

\(\lf\) also differs from the systems of Henkin and Church in its assumptions about infinity. \(\lf\) includes only a rule of \textit{potential} infinity, which says that, for each finite number \(n\), the proposition that there are at least \(n\) things is not contradictory. A rule of \textit{actual} infinity would by contrast assert that for each finite number \(n\), there are \textit{in fact} at least \(n\) things. Potential Infinity suffices for the provability of the \textit{inequalities} we need to be able to prove in mathematics, which we therefore also need to be able to prove throughout the rest of science, such as \(2+2\neq 5\). 

Church and Henkin include axioms of (actual) infinity in their systems, but only for individuals (type $e$\footnote{Or type $\iota$ in Church's and Henkin's terms.}). For this reason, they cannot prove either the potential or actual infinity of propositions, and hence cannot prove \(2+2\neq 5\) on the interpretation of the numerals suitable for counting propositions, which is required, as we saw above, for elementary reasoning about probability. (Even worse, in Henkin's extensional system can prove \(2+2=5\) when the numerals are so interpreted, but that system, as we saw, is also unsuitable for elementary reasoning about probability for an even more basic reason.)

\(\lf\) is closely related to earlier intensional systems, most notably the system \(\mathsf{GM} = \mathsf{ML}_\mathsf{P} + \mathsf{C} + \mathsf{EC} + \mathsf{AC} + \mathsf{PE}\), which is the strongest system that can be assembled out of the components studied by \textcite{Gallin1975-GALIAH} and which is gestured at but not axiomatized by Montague (a similar idea in a more baroque form is presented by \textcite{church1951formulation} as ``Alternative (2)''); a weakening of $\lf$ (with the rule of Potential Infinity restricted to type $e$) is a conservative extension of \(\mathsf{GM}\). Like $\mathsf{GM}$, \(\lf\) significantly improves on the work by Montague (e.g., \cite*{Montague1970-MONUG}) that inspired Gallin's work in its conceptual clarity and economy of notation by eliminating the additional base type of `indices'. $\lf$ is also considerably simpler than $\mathsf{GM}$, in implementing intensionality with the simple rule of Intensionality (attributed, in essence, to \cite{Carnap1942-CARITS-7}: 92 in \cite{Church1943-ReviewOfCarnap}: 300), which permits one to extend a proof of the material equivalence of propositions to a proof of their identity. Other intensional systems that avoid overcomplication by index types, such as Bacon and Dorr's \parencite*{BaconForthcoming-BACC-8} \textit{Classicism}, are too weak in important respects that will be described in future work.

\section{The formal system}

\subsection{Terms}

We begin by introducing the language of $\lf$. It is made up of \textit{terms} that are classified by their \textit{types}. Both the terms and the types by which they are classified are strings (of typographic characters), but for each type $\sigma$ there is also the class of terms of that type, to which the string $\sigma$ will also be used to refer---a tolerable ambiguity that would be eliminated in a fully formalized presentation of the system.

\subsubsection{Types}

There are two \textit{base types}: $e$ (the type of \textit{individual terms}) and $t$ (the type of \textit{sentences} and other \textit{formulae}). The remaining types are \textit{functional types}: for any types $\sigma,\tau$, $\type{\sigma}{\tau}$ is a functional type (which may casually be referred to as `the type of functions from $\sigma$ to $\tau$').

\begin{align}
    \begin{prooftree}
        \hypo{\phantom{\sigma\in\stype}}
        \infer1{e\in\stype}
    \end{prooftree}
    \quad
    \begin{prooftree}
        \hypo{\phantom{\sigma\in\stype}}
        \infer1{t\in\stype}
    \end{prooftree}
    \quad
    \begin{prooftree}
        \hypo{\sigma\in\stype}
        \hypo{\tau\in\stype}
        \infer2{\langle \sigma\tau\rangle \in\stype}
    \end{prooftree}
\end{align}

\subsubsection{Variables}

The terms of the language include, for each type $\sigma$, an infinite stock of variables of that type, each of which is a string consisting of an (italic Roman) letter with the type $\sigma$ as a subscript followed by zero or more primes:
\begin{align}
    \begin{prooftree}
        \hypo{\sigma\in\stype}
        \hypo{\mv{a}\in\mathrm{Letter}}
        \infer2{\mv{a}^\sigma \in\variable}
    \end{prooftree}
    \quad
    \begin{prooftree}
        \hypo{\sigma\in\stype}
        \hypo{\mv{a}\in\mathrm{Letter}}
        \infer2{\mv{a}^\sigma :\sigma}
    \end{prooftree}
    \quad
    \begin{prooftree}
        \hypo{\mv{a}\in\variable}
        \infer1{\mv{a}'\in\variable}
    \end{prooftree}
    \quad
    \begin{prooftree}
        \hypo{\mv{a}\in\variable}
        \hypo{\mv{a}:\sigma}
        \infer2{\mv{a}':\sigma}
    \end{prooftree}
\end{align}

\subsubsection{Constants}

The \textit{constants} of the language are the strings consisting of the standard inclusion symbol $\subseteq$ followed by a type $\sigma$ as a subscript, which indicates the type of the constant according to the rule that the type of $\subseteq_\sigma$ is $\langle\sigma t\rangle\langle \sigma t \rangle t$:
\begin{align}
    \begin{prooftree}
        \hypo{\sigma\in\stype}
        \infer1{\subseteq_\sigma:\langle\sigma t\rangle\langle \sigma t \rangle t}
    \end{prooftree}
\end{align}
\(\subseteq_\sigma\) is pronounced ``every'', or, when written in infix notation (\cref{subsec:infix}), ``implies'' or ``is included in'' (and variants).

\subsubsection{Complex terms}

All \textit{complex terms} of the language are constructed in accordance with one of two formation rules: (function) application and (function) abstraction:

\begin{align}
    \begin{prooftree}
        \hypo{\mv{f}:\sigma\tau}
        \hypo{\mv{a}:\sigma}
        \infer2[\scriptsize{App}]{(\mv{fa}):\tau}
    \end{prooftree}
    \quad\quad
    \begin{prooftree}
        \hypo{\mv{x}\in\variable}
        \hypo{\mv{x}:\sigma}
        \hypo{\mv{A}:\tau}
        \infer3[\scriptsize{Abs}]{(\lambda \mv{x}\ldot \mv{A}):\sigma \tau}
    \end{prooftree}
\end{align}

\subsubsection{$\beta$-equivalence}\label{Sec-BetaEquivDef}

Terms of the forms $(\lambda \mv{x}\ldot\mv{A})\mv{B}$ and $[\mv{B}/\mv{x}]\mv{A}$ are said to be \textit{immediately $\beta$-equivalent}, and terms are said to be \textit{$\beta$-equivalent} when one can be obtained from the other by any number of substitutions of immediately $\beta$-equivalent parts (each term counts as part of itself). We write $\mv{A}\sim_\beta \mv{B}$ when $\mv{A}$ and $\mv{B}$ are $\beta$-equivalent. 

\subsubsection{Expressions}

An \textit{expression} of the language is a term with no free variables.

\subsubsection{Formulae and sentences}

A \textit{formula} is a term of type \(t\). A \textit{sentence} is an expression of type \(t\).

\subsection{Abbreviations and notational conventions}\label{sec:abbreviations}

\subsubsection{Omission of parentheses and type decorations}

Parentheses and angle brackets ($\langle,\rangle$) may be omitted in accordance with the following conventions.
\begin{enumerate}
    \item Variable binders (i.e., \(\lambda\) and other introduced variable binders; which are indicated by a \(\ldot\) following the bound variable) take widest possible scope.
    \item Boolean connectives ($\wedge,\vee,\rightarrow,\leftrightarrow$) are written in infix notation and take next widest scope.
    \item Other function terms written in infix notation take next greatest scope.
    \item Functions written in prefix notation take smallest possible scope.
    \item Functions written in infix notation associate to the right, when the conventions above do not say otherwise (we apply this convention only when it is clear that the functions in question provably associate). So, given all of the foregoing,
\begin{align}
    \lambda p\ldot fp\rightarrow q
\end{align}
omits parentheses which should be filled in as
\begin{align}
    (\lambda p\ldot ((fp)\rightarrow q)).
\end{align}
    \item Angle brackets (\(\langle,\rangle\)) are omitted in accordance with the convention that complex types associate to the right. So, e.g., \(ttt\) abbreviates \(\langle t\langle tt \rangle\rangle\).
\end{enumerate}

Type decorations may be omitted when writing terms when they can be inferred from context, i.e., when there is a unique term that could result in the written string by deleting type decorations from any variables or constants, subject to the following conventions.
\begin{enumerate}
    \item Letters that could be bound by the same variable binder must have the same type. E.g., \(\lambda x^e\ldot x\) and \(\lambda x\ldot x^e\) must be \(\lambda x^e\ldot x^e\).
    \item \(p\) and \(q\) always have type \(t\) unless otherwise decorated.
\end{enumerate}
For example, 
\begin{align}
    \lambda x^e\ldot ffx
\end{align}
abbreviates
\begin{align}
    \lambda x^e\ldot f^{ee}f^{ee}x^e.
\end{align}

\subsubsection{Infix Notation}\label{subsec:infix}

Introduced constants \(\mv{f}\) of types that take at least two arguments, i.e., types of the form \(\sigma\tau\rho\), may be written in the form
\begin{align}
    (\mv{a}\;\mv{f}\;\mv{b})\rewrite ((\mv{f}\mv{a})\mv{b})
\end{align}

\subsubsection{Metanotation}\label{subsec:meta}

Bold symbols are used as metavariables. \(\mv{P}\) and \(\mv{Q}\) always range over formulae. The Greek letters \(\sigma\), \(\tau\), and \(\rho\) are also metavariables, but always range over types.

The overset arrow \(\vec{\cdot}\) is used atop a metavariable to construct a metavariable ranging over strings of variables or types, of arbitrary length. When this device is used to range over both a string of variables and a string of types, as in
\begin{align}
    \lambda X^{\vec{\sigma}t}\vec{z}\ldot \neg X\vec{z},
\end{align}
the intended instances are those in which the variables \(\vec{z}\) have the respective types \(\vec{\sigma}\).

\subsubsection{Basic definitions}
\begin{align}
    \top &\rewrite (\lambda p\ldot p)\subseteq (\lambda p\ldot p)\\
    \forall_\sigma &\rewrite \lambda X^{\sigma t}\ldot ((\lambda y^\sigma \ldot \top)\subseteq X)\\
    \forall \mv{x}^\sigma \ldot\mv{P}&\rewrite \forall_\sigma \lambda \mv{x}^\sigma\ldot \mv{P}\\
    \bot &\rewrite (\lambda p\ldot \top)\subseteq (\lambda p\ldot p)\\
    \neg &\rewrite \lambda p\ldot ((\lambda r^t\ldot p)\subseteq (\lambda r^t\ldot \bot))\\
    \exists_\sigma &\rewrite \lambda X^{\sigma t}\ldot \neg \forall y^\sigma \ldot \neg Xy\\
    \exists \mv{x}^\sigma \ldot\mv{P}&\rewrite \exists_\sigma \lambda \mv{x}^\sigma\ldot \mv{P}\\
    \lambda \mv{x}\vec{\mv{y}}\ldot\mv{A}&\rewrite \lambda \mv{x}\ldot\lambda \vec{\mv{y}}\ldot \mv{A}\\
    \forall \mv{x}\vec{\mv{y}} \ldot\mv{P}&\rewrite \forall \mv{x}\ldot \forall \vec{\mv{y}}\ldot \mv{P}\\
    \exists \mv{x}\vec{\mv{y}} \ldot\mv{P}&\rewrite \exists \mv{x}\ldot \forall \vec{\mv{y}}\ldot \mv{P}\\
    \rightarrow &\rewrite \lambda p q\ldot ((\lambda r^t\ldot p)\subseteq (\lambda r^t\ldot q))\label{eq:defrightarrow}\\
    \vee &\rewrite \lambda pq\ldot (\neg p\rightarrow q)\\
    \wedge &\rewrite \lambda pq\ldot \neg (\neg p\vee \neg q)\\
    \leftrightarrow &\rewrite \lambda pq\ldot ((p\rightarrow q)\wedge (q\rightarrow p))\\
    \equiv_{\vec\sigma t}&\rewrite \lambda XY^{\vec\sigma t}\vec{z}\ldot (X\vec{z}\leftrightarrow Y\vec{z})\\
    =_\sigma&\rewrite \lambda xy^\sigma\ldot ((\lambda Z\ldot Zx)\subseteq (\lambda Z\ldot Zy))
    \\
    \neq_\sigma&\rewrite\lambda xy\ldot \neg (x=_\sigma y)
    \\
    \nec &\rewrite {=}\top
    \\
    \pos &\rewrite {\neq}\bot
    \\
    0_\sigma &\rewrite \lambda X^{\sigma t}\ldot \neg \exists X\label{eq:def0}\\
    1_\sigma &\rewrite \lambda X^{\sigma t}\ldot \exists y\ldot (Xy\wedge \forall z\ldot (Xz\rightarrow y=z))\label{eq:def1}\\
    \setminus_{\vec{\sigma} t}&\rewrite \lambda XY^{\vec{\sigma}t}\vec{z}\ldot (X\vec{z}\wedge \neg Y\vec{z})\label{eq:setminus}
    \\
    +_\sigma &\rewrite \lambda mn^{\langle \sigma t\rangle t}X^{\sigma t}\ldot \exists Y\ldot (Y\subseteq X\wedge mY\wedge n(X\setminus Y))\label{eq:def+}\\
    \nat_\sigma &\rewrite \lambda n^{\langle \sigma t\rangle t}\forall X\ldot (X0\rightarrow (X\subseteq \lambda y\ldot X(y+1)) \rightarrow Xn)
\end{align}

\subsection{\(\lf\)}

\(\lf\) is presented as a natural deduction system in sequent calculus notation. A \textit{sequent} is a string consisting of a possibly empty list of formulae (separated by commas), followed by the symbol \(\proves\), and then a formula. The sequent
\begin{align}
    \Gamma\proves \mv{P}
\end{align}
may be read ``\(\mv{P}\) has been proved (or derived) from the assumptions \(\Gamma\)''.

A \textit{rule} is written in the following form, where capital Greek letters range over lists of formulae (i.e., formulae separated by commas),
\begin{align}
    \begin{prooftree}
        \hypo{\Gamma_1\proves \mv{P_1}}
        \hypo{\dots}
        \hypo{\Gamma_n\proves \mv{P_n}}
        \infer3{\Delta\proves \mv{Q}}
    \end{prooftree}
\end{align}
and it indicates that, when \(\mv{P_i}\) is provable from \(\Gamma_i\) for each \(i\in \{1, \dots, n\}\), \(\mv{Q}\) is provable from \(\Delta\).

\(\lf\) is defined as a class of rules, namely the rules \labelcref{rule:structural,rule:beta,rule:universalinstantiation,rule:universalgeneralization,rule:Intensionality,rule:functionextensionality,rule:Choice,rule:potentialinfinity}. A formula \(\mv{P}\) is \textit{provable} in \(\lf\) from the assumptions \(\Gamma\) if and only if the sequent
\begin{align}
    \Gamma\proves\mv{P}
\end{align}
is derivable by some application of these rules. A \textit{theorem} of \(\lf\) is a sentence that is provable from no assumptions (i.e., from the empty list of assumptions, in which case we write: $\proves\mv{P}$).

\renewcommand{\thesubsubsection}{R.\arabic{subsubsection}}

\subsubsection{Structural Rules}\label{rule:structural}

\begin{align}
    \begin{prooftree}
        \hypo{\vphantom{\Gamma,\mv{P},\mv{P}\proves\mv{Q} }}
        \infer1{\Gamma,\mv{P}\proves\mv{P} }
    \end{prooftree}
    \quad\quad
    \begin{prooftree}
        \hypo{\Gamma,\mv{P},\mv{P}\proves\mv{Q} }
        \infer1{\Gamma,\mv{P}\proves\mv{Q} }
    \end{prooftree}
    \quad\quad
    \begin{prooftree}
        \hypo{\Gamma\proves\mv{Q} }
        \infer1{\Gamma,\mv{P}\proves\mv{Q} }
    \end{prooftree}
    \quad\quad
    \begin{prooftree}
        \hypo{\Gamma,\mv{P},\mv{Q},\Delta\proves\mv{R}}
        \infer1{\Gamma,\mv{Q},\mv{P},\Delta\proves\mv{R}}
    \end{prooftree}
    \quad\quad
    \begin{prooftree}
        \hypo{\Gamma\proves \mv{P}}
        \hypo{\Delta,\mv{P}\proves \mv{Q}}
        \infer2{\Gamma,\Delta\proves\mv{Q}}
    \end{prooftree}
\end{align}

\subsubsection{\(\beta\)-Equivalence}\label{rule:beta}

\begin{align}
    \begin{prooftree}
        \hypo{\Gamma\proves \mv{P}}
        \infer1{\Gamma\proves\mv{Q}}
    \end{prooftree}
\end{align}
where \(\mv{P}\sim_\beta \mv{Q}\) (see \S\ref{Sec-BetaEquivDef}). 

\subsubsection{Universal Instantiation}\label{rule:universalinstantiation}

\begin{align}
    \begin{prooftree}
        \hypo{\Gamma\proves \mv{F}\subseteq_\sigma \mv{G}}
        \hypo{\Gamma\proves \mv{Fa}}
        \infer2{\Gamma\proves \mv{Fb}}
    \end{prooftree}
\end{align}

\subsubsection{Universal Generalization}\label{rule:universalgeneralization}

\begin{align}
    \begin{prooftree}
        \hypo{\Gamma,\mv{Fx}\proves \mv{Gx}}
        \infer1{\Gamma\proves \mv{F}\subseteq_\sigma \mv{G}}
    \end{prooftree}
\end{align}
where \(\mv{x}\) is a variable of type \(\sigma\) not free in $\mv{F}$, $\mv{G}$, or any formula in \(\Gamma\).

\subsubsection{Negation Elimination}\label{rule:negationelimination}
\begin{align}
    \begin{prooftree}
        \hypo{\Gamma, \neg\mv{P}\proves \mv{P}}
        \infer1{\Gamma\proves \mv{P}}
    \end{prooftree}
\end{align}



\subsubsection{Intensionality}\label{rule:Intensionality}

\begin{align}
    \begin{prooftree}
        \hypo{ \mv{P}\proves \mv{Q}}
        \hypo{ \mv{Q}\proves \mv{P}}
        \infer2{\proves {\mv{P}=\mv{Q}}}
    \end{prooftree}
\end{align}


\subsubsection{Function Extensionality}\label{rule:functionextensionality}
\begin{align}
    \begin{prooftree}
        \hypo{ \Gamma\proves \mv{f}\mv{x}=_\tau \mv{g}\mv{x}}
        \infer1{ \Gamma\proves {\mv{f}=_{\sigma\tau}\mv{g}}}
    \end{prooftree}
\end{align}
where \(\mv{x}\) is a variable of type \(\sigma\) not free in \(\mv{f}\), \(\mv{g}\), or any formula in \(\Gamma\).

\subsubsection{Choice}\label{rule:Choice}

\begin{align}
    \begin{prooftree}
        \hypo{\Gamma \proves \forall x^\sigma \ldot\exists y^\tau \ldot\mv{R}xy}
        \infer1{\Gamma \proves \exists f^{\sigma \tau}\ldot\forall x\ldot\mv{R}x(fx)}
    \end{prooftree}
\end{align}

\subsubsection{Potential Infinity}\label{rule:potentialinfinity}
\begin{align}
    \begin{prooftree}
        \hypo{\Gamma \proves \nat_\sigma \mv{n}}
        \infer1{\Gamma \proves \bot\neq \exists_{\sigma t} \mv{n}}
    \end{prooftree}
\end{align}
where \(\sigma\) is either \(e\) or \(t\).

\renewcommand{\thesubsubsection}{\thesection.\thesubsection.\arabic{subsubsection}}

\section{Adding and subtracting assumptions (\({}+X\); \({}-\mathrm{R.n}\))}

By ``\(\lf-\mathrm{R.n}\)'' we mean the system that has the rules of \(\lf\) excluding R.n (where n is among 1 through 9).

Where \(X\) is a class of formulae, \(\lf+X\), is the system whose proofs are proofs in \(\lf\) from assumptions in \(X\). That is, a formula \(\mv{P}\) is a theorem of \(\lf+X\) if and only if there is a proof in \(\lf\) that concludes with the sequent
\begin{align}
    \Gamma\proves \mv{P},
\end{align}
where \(\Gamma\) is a list of formulae in \(X\). Note that \(\lf+X\) is not necessarily closed under the rule of Intensionality, i.e., \(\lf+X\) may prove the extensional equivalence \(\mv{P}\leftrightarrow \mv{Q}\) but not the identity \(\mv{P}=\mv{Q}\).

\section{\(\lf_\iota\); \(\lf_\varepsilon\)}\label{Sec-IotaAndEpsilon}

We introduce two new families of constants indexed by simple types \(\sigma\):
\begin{align}
    \iota_\sigma &: \langle \sigma t\rangle\sigma\\
    \varepsilon_\sigma &: \langle \sigma t\rangle \sigma 
\end{align}
\(\lf_\iota\) and \(\lf_{\varepsilon}\) are used for reasoning with these constants.

\subsection{\(\lf_\iota\)}

\(\iota\) is Church's (\cite{Church1940-CHUAFO}: 61) \textit{description function}: a function which maps every uniquely instantiated property to an instance thereof. \(\mathrm{D}_\iota\) is a set of axioms to this effect, which also choose a ``default'' value for \(\iota\) to take, \(\dagger\), when applied to properties which are not uniquely instantiated:
\begin{gather}
    \forall X^{\sigma t}\ldot (\exists! X\rightarrow X(\iota X))\label{eq:desc1}\\
    \forall X^{\sigma t}\ldot (\neg\exists !X\rightarrow \iota X=\dagger)\label{eq:desc2}
\end{gather}
\(\dagger\) is defined by the following recursive rewrite rules (although the particular choice of \(\dagger\) is immaterial):
\begin{align}
    \dagger_e&\rewrite \iota \lambda x^e\ldot \bot\\
    \dagger_t&\rewrite \bot\\
    \dagger_{\sigma\tau}&\rewrite \lambda x^{\sigma}\ldot \dagger_\tau
\end{align}
\(\lf_\iota\) is \(\lf+\mathrm{D}_\iota\), i.e., \(\lf\) plus every sentence of the form \labelcref{eq:desc1} or \labelcref{eq:desc2}.

The utility of \(\lf_\iota\) should be obvious: the language of mathematics is replete with notations (especially, but not only, variable-binding operators\footnote{See \cite{KalishMontague1964-LTFL}: Ch. IX for a representative collection of examples.}) that are definable in \(\lf_\iota\) but not in \(\lf\). Perhaps the most common example is class abstraction, which is definable by
\begin{align}
    \{x^\sigma :\mv{P}\}\rewrite \iota \lambda X^{\sigma t}\ldot(\forall y\ldot(Xy=\top\vee Xy=\bot)\wedge X\equiv\lambda x\ldot \mv{P}).\label{Def-ClassAbs}
\end{align}

\subsection{\(\lf_{\varepsilon}\)}

\(\varepsilon\) is Hilbert's (\cite*{Hilbert1922-Abhandlungen}) \textit{epsilon operator}\footnote{\textcite{Hilbert1922-Abhandlungen} does not use the lowercase epsilon ($\varepsilon$) exactly our way. Rather, he uses it as a variable-binding operator (the formation rule for such an operator would be: if $\mv{x}$ is a variable of type $\sigma$, then $(\varepsilon\mv{x}\ldot \mv{P})$ is a term of type $\sigma$). In the present setting, such a variable-binding operator would be redundant: were it to be to be introduced, the $\varepsilon$-function could be defined in terms of it by
\begin{align}
    \varepsilon_\sigma \rewrite \lambda X^{\sigma t}\ldot\varepsilon y^\sigma \ldot Xy,
\end{align}
while, in following Church (\cite*{Church1940-CHUAFO}: 57, 61, where $\iota$ is confusingly used for both the description and Choice functions) in adopting the $\varepsilon$-operators as constants, we can (and we do) also follow Church (\cite*{Church1940-CHUAFO} : 58) in eliminating the corresponding variable-binding operator, as we do others, by
\begin{align}
     (\varepsilon x^\sigma \ldot \mv{P}) \rewrite \varepsilon_\sigma \lambda x^\sigma \ldot\mv{P}.
\end{align}} and Church's (\cite*{Church1940-CHUAFO}: 61) \textit{Choice function}. It maps every instantiated property (whether or not uniquely instantiated) to an instance thereof, and everything else to \(\dagger\). This is captured by the class \(\mathrm{C}_\varepsilon\), which consists of the following sentences, where \(\sigma\) is any type.
\begin{gather}
    \forall X^{\sigma t}\ldot (\exists X\rightarrow X(\varepsilon X))\label{eq:eps1}\\
    \forall X^{\sigma t}\ldot (\neg\exists X\rightarrow \varepsilon X=\dagger)\label{eq:eps2}
\end{gather}
\(\lf_\varepsilon\) is \(\lf+\mathrm{D}_\iota+\mathrm{C}_\varepsilon\).\footnote{It is redundant to have the primitive description functions \(\iota\), since \(\lf_\varepsilon\) categorically reduces to \(\lf+\mathrm{C}_\varepsilon\) by the definition schema
\begin{align}
    \dagger_e&\rewrite \varepsilon \lambda x^e\ldot \bot\\
    \iota_\sigma &\rewrite \varepsilon\lambda f^{\langle \sigma t\rangle \sigma}\ldot \forall X^{\sigma t}\ldot ((\exists ! X\rightarrow X(fX))\wedge (\neg \exists!X\rightarrow fX=\dagger_\sigma))
\end{align}
However, the use of \(\iota\) is convenient to conceptually distinguish between full-blown applications of Choice and mere applications of descriptions.}

The utilty of \(\lf_\varepsilon\) is less obvious than the utility of \(\lf_\iota\), and, in fact, it is hard to deny that the deductive power \(\lf_\varepsilon\) adds to that of \(\lf_\iota\) is less often useful than the deductive power that \(\lf_\iota\) adds to that of \(\lf\). However, the former deductive power is sometimes useful, and that is enough to justify the adoption of \(\lf_\varepsilon\). As Carnap (\cite*{Carnap1962-UHEOST}) points out, $\varepsilon$s can be used to extract definitions from non-categorical theories, $\varepsilon\mv{F}$ being a term that defines a unique satisfier of $\mv{F}$ provided that there is at least one satisfier of $\mv{F}$. As a concrete example of where this might be useful, in future work on the foundations of semantics, we show that various useful semantic (or proto-semantic) notions are definable in \(\lf_\varepsilon\) + syntax but not in \(\lf_\iota\) + syntax.

\subsection{The slingshot argument}\label{Sec-Slingshot}

We stress that \(\lf_\iota\) and \(\lf_{\varepsilon}\) are \textit{not} closed under the rule of Intensionality (if \(\lf\) is consistent), and consequently there are formulae \(\mv{P}\) and \(\mv{Q}\) that are provably equivalent in \(\lf_\iota\) (hence also \(\lf_{\varepsilon}\)), but for which
\begin{align}
    \mv{P}=\mv{Q}
\end{align}
is not provable, and indeed is refutable,\footnote{Let\begin{align}
    \alpha\rewrite\forall p\ldot(p\bi@p),
\end{align}with @ defined as in (\ref{Def-Actuality}). Then $\alpha\bi\top$ is provable but $\alpha=\top$ is refutable in these systems.}
in the same system.

For example, in \(\lf_\iota\), consider the function \(@\), defined as
\begin{align}
    @&\rewrite \lambda p\ldot \iota q\ldot ((p\rightarrow q=\top) \wedge (\neg p\rightarrow q=\bot)). \label{Def-Actuality}
\end{align}
\(@\) is the function that maps every truth to \(\top\) and every falsehood to \(\bot\), so the material equivalence of \(@p\) with \(p\) is intuitively clear (it is also easily provable). However, the provability of the identity
\begin{align}
    p=@p\label{Fmla-Collapse}
\end{align}
would have a collapsing effect: since \(@\) only ever takes the values \(\top\) or \(\bot\), the provability of (\ref{Fmla-Collapse}) would result in the provability of \(\top\) and \(\bot\) being the only propositions. This is in contradiction with Potential Infinity (\labelcref{rule:potentialinfinity}; see \cref{metathm:nonextensionality}).

Inferring the identity of \(p=@p\) (or some close variant\footnote{Specifically, variants constructed using directly by means of $\iota$-terms (such as $\iota \lambda q\ldot((p\to q=\top)\wedge (p\to q=\top))$) or by means of class abstracts (such as $\{ q:(p\to q=\top)\wedge (p\to q=\bot)\}$), which are definable by means of $\iota$ (e.g., by (\ref{Def-ClassAbs}) or, equivalently, with $\{x:\mv{P}\}$ taken to abbreviate $\lambda x\ldot @\mv{P}$). In the literature, these tend to be preferred to $@p$.\label{Fn-SlingsgotVariants}}) from the proof of the material equivalence \(p\bi @p\) is known as the \textit{slingshot} argument. It is deployed, among others, by Church (\cite*{Church1943-ReviewOfCarnap}: 299--301) in his review of Carnap's (\cite*{Carnap1942-CARITS-7}) \textit{Introduction to Semantics}, where Church attributes it to \textcite{Frege1892-SinnBedeutung}; by Gödel (\cite*{Godel1944-RsML}: 450, esp. n. 5) in his ``Russell's Mathematical Logic''; again by Church (\cite*{Church1956-IML}: 24--25) in his \textit{Introduction to Mathematical Logic}, where it is again attributed to \textcite{Frege1892-SinnBedeutung}; by \textcite{Quine1953-RefAndMod} in his criticism of quantified modal logic; and by Myhill (\cite*{Myhill1958-MYHPAI}: 77--78) in relation to Church's \parencite*{church1951formulation} version of the rule of Intensionality (called \textit{Alternative (2)}). The slingshot is not valid in \(\lf_\iota\) or \(\lf_{\varepsilon}\), so these systems avoid the collapse.

\section{Basic facts}

\begin{metatheorem}[Modus Ponens and Conditional Proof]
    \(\lf\) is closed under the rules of Modus Ponens and Conditional Proof.
    \begin{align}
        \begin{prooftree}
            \hypo{\Gamma\proves \mv{P}}
            \hypo{\Gamma\proves\mv{P}\rightarrow\mv{Q}}
            \infer2{\Gamma\proves \mv{Q}}
        \end{prooftree}
        \quad\quad
        \begin{prooftree}
            \hypo{\Gamma,\mv{P}\proves \mv{Q}}
            \infer1{\Gamma\proves \mv{P}\rightarrow\mv{Q}}
        \end{prooftree}
    \end{align}
\end{metatheorem}
\begin{proof}
    Assume without loss of generality that \(p,q,r\) do not occur in \(\Gamma,\mv{P},\mv{Q}\) (otherwise use different variables). For Modus Ponens we have
    \begin{align}
        \begin{prooftree}
            \hypo{\Gamma\proves \mv{P}}
            \infer1[\scriptsize \(\beta\)]{\Gamma\proves (\lambda r\ldot \mv{P})p}
            \hypo{\Gamma\proves \mv{P}\rightarrow\mv{Q}}
            \infer1[\scriptsize Def.\ \labelcref{eq:defrightarrow}]{\Gamma\proves (\lambda r\ldot \mv{P})\subseteq(\lambda r\ldot \mv{Q})}
            \infer2[\scriptsize Universal Instantiation]{\Gamma\proves(\lambda r\ldot \mv{Q})p}
            \infer1[\scriptsize \(\beta\)]{\Gamma\proves\mv{Q}}
        \end{prooftree}
    \end{align}
    for Conditional Proof we have
    \begin{align}
        \begin{prooftree}
            \hypo{\Gamma,\mv{P}\proves\mv{Q}}
            \infer0{\Gamma,(\lambda r\ldot \mv{P})p\proves (\lambda r\ldot \mv{P})p}
            \infer1[\scriptsize \(\beta\)]{\Gamma,(\lambda r\ldot \mv{P})p\proves \mv{P}}
            \infer2{\Gamma,(\lambda r\ldot \mv{P})p\proves \mv{Q}}
            \infer1[\scriptsize \(\beta\)]{\Gamma,(\lambda r\ldot \mv{P})p\proves (\lambda r\ldot \mv{Q})p}
            \infer1[\scriptsize Universal Generalization]{\Gamma\proves (\lambda r\ldot \mv{P})\subseteq (\lambda r\ldot \mv{Q})}
            \infer1[\scriptsize Def.\ \labelcref{eq:defrightarrow}]{\Gamma\proves \mv{P}\rightarrow\mv{Q}}
        \end{prooftree}
    \end{align}
\end{proof}

\begin{metatheorem}[Intuitionistic Propositional Logic]\label{metathm:intuitionistic}
    Every theorem of intuitionistic propositional logic is provable from rules \labelcref{rule:structural,rule:beta,rule:universalinstantiation,rule:universalgeneralization}.
\end{metatheorem}

\begin{metatheorem}[Classical Quantification Theory]\label{metathm:classicalquantified}
    Classical introduction and elimination rules for \(\forall_\sigma\) are derivable from
\labelcref{rule:structural,rule:beta,rule:universalinstantiation,rule:universalgeneralization}.
\end{metatheorem}

\begin{metatheorem}[The Logic of Identity]\label{metathm:identity}
    The reflexivity of identity and Leibniz' law are provable from \labelcref{rule:structural,rule:beta,rule:universalinstantiation,rule:universalgeneralization}.
    \begin{gather}
        \forall x^\sigma\ldot x= x
        \\
        \forall xy^\sigma \ldot (x= y\rightarrow \mv{P}\rightarrow [y/x]\mv{P})
    \end{gather}
\end{metatheorem}
\begin{proof}
    For Leibniz' law, we first unpack the definition of \(=\) in the antecedent as antecedent to get (with \labelcref{rule:beta})
    \begin{align}
        (\lambda Z\ldot Zx)\subseteq (\lambda Z\ldot Zy)
    \end{align}
    which, by universal instantiation (\labelcref{rule:universalinstantiation}) in conjunction with
    \begin{align}
        (\lambda x\ldot \mv{P})x,
    \end{align}
    yields
    \begin{align}
        (\lambda x\ldot \mv{P})y,
    \end{align}
    from which Leibniz' law follows by conditional proof and Universal Generalization.
\end{proof}

\begin{metatheorem}[Intensionality]\label{metathm:fullIntensionality}
    \(\lf\) is intensional, i.e., when \(\mv{P}\) and \(\mv{Q}\) are provable from each other, \(\mv{P}\) and \(\mv{Q}\) can be substituted in any context, including cases with variable capture.
\end{metatheorem}
\begin{proof}
    Let \(\vec{\mv{x}}\) be the variables that occur free in \(\mv{P}\) or \(\mv{Q}\). From Intensionality and Function Extensionality we have
    \begin{align}
        \lambda \vec{\mv{x}}\ldot \mv{P}=\lambda \vec{\mv{x}}\ldot \mv{Q}
    \end{align}
    which are closed terms, so can be substituted in any context by Leibniz' law (\cref{metathm:identity}). And by \(\beta\)-Equivalence (\labelcref{rule:beta}), \(\lambda \vec{\mv{x}}\ldot \mv{P}\) can be substituted in any context with \((\lambda \vec{\mv{x}}\ldot \mv{P})\vec{\mv{x}}\), and similarly for \(\mv{Q}\).
\end{proof}

\begin{metatheorem}[Classical Propositional Logic]\label{metathm:classicalpropositional}
    Every theorem of classical propositional logic is provable from rules \labelcref{rule:structural,rule:beta,rule:universalinstantiation,rule:universalgeneralization,rule:negationelimination}.
\end{metatheorem}
\begin{proof}
    \labelcref{rule:negationelimination} is sufficient to derive classical propositional logic from intuitionistic propositional logic, which is in turn supplied by
    \cref{metathm:intuitionistic}.
\end{proof}


\begin{metatheorem}[The Modal Logic $\mathsf{S4}$]\label{metathm:S4}
    Every theorem of the propositional modal logic \(\mathsf{S4}\) is provable from the rules \labelcref{rule:structural,rule:beta,rule:universalinstantiation,rule:universalgeneralization,rule:negationelimination,rule:Intensionality,rule:functionextensionality}. Moreover, the rule of necessitation,
    \begin{align}
        \begin{prooftree}
            \hypo{\proves\mv{P}}
            \infer1{\proves\nec\mv{P}}
        \end{prooftree}
    \end{align}
    is derivable from those rules.
\end{metatheorem}

\begin{metatheorem}[The Modal Logic $\mathsf{S5}$ ]\label{metathm:S5}
    The ``5'' axiom of modal logic, 
\begin{align}
    \forall p\ldot (\neg\nec p\rightarrow \nec\neg \nec p),
\end{align}
is provable from rules \labelcref{rule:structural,rule:beta,rule:universalinstantiation,rule:universalgeneralization,rule:negationelimination,rule:Intensionality,rule:functionextensionality,rule:Choice}. Thus, \(\lf\) includes the modal logic \(\mathsf{S5}\).
\end{metatheorem}
\begin{proof}[Proof of \cref{metathm:S4,metathm:S5}]
    \(\nec p\) is by definition \(p=\top\), so we must show
    \begin{align}
        p\neq \top\rightarrow \nec(p\neq \top).
    \end{align}
    We have for any function \(f\),
    \begin{align}
        \nec(fp\neq f\top\rightarrow p\neq \top).
    \end{align}
    Hence (by basic modal logic; see \cref{metathm:S4})
    \begin{align}
        \nec(fp\neq f\top)\rightarrow \nec(p\neq \top).
    \end{align}
    So it suffices to find some function \(f\) for which
    \begin{align}
        p\neq \top\rightarrow \nec(fp\neq f\top),
    \end{align}
    i.e., a function that, for all \(p\) distinct from \(\top\), maps \(p\) and \(\top\) to any two necessarily distinct things, such as, \(0\) and \(1\) or \(\bot\) and \(\top\). This is supplied by Choice applied to the relation
    \begin{align}
        \lambda pq\ldot ((p\neq \top\rightarrow q=\bot)\wedge(p=\top\rightarrow q=\top)).
    \end{align}
\end{proof}

\begin{metatheorem}[Necessity of Identity and Distinctness]
Every instance of the necessity of identity,
    \begin{align}
        \forall xy^\sigma\ldot (x=y\rightarrow \nec(x=y)),
    \end{align}
    is provable from rules \labelcref{rule:structural,rule:beta,rule:universalinstantiation,rule:universalgeneralization,rule:Intensionality}, and every instance of the necessity of distinctness, 
    \begin{align}
        \forall xy^\sigma\ldot (x\neq y\rightarrow \nec (x\neq y)),
    \end{align}
    is provable from rules \labelcref{rule:structural,rule:beta,rule:universalinstantiation,rule:universalgeneralization,rule:negationelimination,rule:Intensionality,rule:Choice}.
\end{metatheorem}
\begin{proof}
    For necessity of identity, we have \(\nec (x=x)\) by intensionality and hence necessity of identity by Leibniz' law (\cref{metathm:identity}). For necessity of distinctness, we use Choice as in \cref{metathm:S5}.
\end{proof}

\begin{metatheorem}[The Barcan Formula and its Converse; \cite{Bacon2018-BACTBN}]\label{metathm:barcan}
    Every instance of the Barcan formula and its converse is provable from rules \labelcref{rule:structural,rule:beta,rule:universalinstantiation,rule:universalgeneralization,rule:negationelimination,rule:Intensionality,rule:functionextensionality}.
    \begin{align}
        \forall X^{\vec\sigma t} \ldot (\forall \vec{z}\ldot \nec X\vec{z})\rightarrow \nec \forall \vec{z}\ldot X\vec{z}\label{eq:bf}\\
        \forall X^{\vec\sigma t}\ldot (\nec \forall \vec{z}\ldot X\vec{z})\rightarrow \forall \vec{z}\ldot \nec X\vec{z}\label{eq:cbf}
    \end{align}
\end{metatheorem}
\begin{proof}
    For \labelcref{eq:bf} (the Barcan formula) in the case where \(\vec{\sigma}\) has length one, we observe 
    \begin{align}
        (\forall z\ldot \nec Xz)\rightarrow \forall z\ldot (Xz=(\lambda y\ldot \top)z) 
    \end{align}
    hence, by \labelcref{rule:functionextensionality},
    \begin{align}
        X=\lambda y\ldot \top.
    \end{align}
    Moreover, by \labelcref{rule:Intensionality} we have
    \begin{align}
        \nec \forall z\ldot (\lambda y\ldot \top)z
    \end{align}
    hence \(\nec\forall z\ldot Xz\) as required. This proof extends to \(\vec{\sigma}\) of arbitrary length by induction.

    For \labelcref{eq:cbf} (the converse Barcan formula), we have by \labelcref{rule:Intensionality} the identity
    \begin{align}
        (\forall \vec{z}\ldot X\vec{z})=(X\vec{y}\wedge \forall \vec{z}\ldot X\vec{z}).
    \end{align}
    we also have \(p= (p\wedge \top) \),
    hence
    \begin{align}
        (\nec \forall \vec{z}\ldot X\vec{z})\rightarrow \nec(X\vec{y})
    \end{align}
    which yields \labelcref{eq:cbf} by basic quantificational reasoning.
\end{proof}

\begin{metatheorem}[Modal Intensionality; \cite{Bacon2018-BACTBN}]\label{metathm:intensionalism}
    All axioms of modal intensionality---asserting that necessarily equivalent propositions and necessarily equivalent properties are identical---are provable from rules \labelcref{rule:structural,rule:beta,rule:universalinstantiation,rule:universalgeneralization,rule:negationelimination,rule:Intensionality,rule:functionextensionality}.
    \begin{gather}
        \forall pq\ldot (\nec(p\leftrightarrow q)\rightarrow p=q)\label{eq:propositionintensionalism}\\
        \forall FG^{\vec{\sigma}t}\ldot (\nec (F\equiv G)\rightarrow F=G)\label{eq:propertyintensionalism}
    \end{gather}
\end{metatheorem}
\begin{proof}
    By some basic applications of \labelcref{rule:Intensionality} we have the equations 
    \begin{gather}
        p=(p\wedge \top)\label{eq:pptop}\\
        (p\wedge q)=(p\wedge (p\rightarrow q))\label{eq:pqparrowq}.
    \end{gather}
    By \labelcref{eq:pptop} we have
    \begin{align}
        \nec(p\leftrightarrow q)\rightarrow (p=(p\wedge (p\rightarrow q))\wedge q=(q\wedge (q\rightarrow p))),
    \end{align}
    hence, by \labelcref{eq:pqparrowq}
    \begin{align}
        \nec (p\leftrightarrow q)\rightarrow (p=(p\wedge q)\wedge q=(p\wedge q)),
    \end{align}
    which implies \labelcref{eq:propositionintensionalism}.
    
    For \labelcref{eq:propertyintensionalism}, we prove the case where \(\vec{\sigma}\) has length one, i.e.,
    \begin{align}
        \forall FG^{\sigma t}\ldot (\nec (F\equiv G)\rightarrow F=G).
    \end{align}
    The other cases follow straightforwardly. 
    For this, suppose \(\nec(F\equiv G)\), i.e.,
    \begin{align}
        \nec \forall x^\sigma \ldot (Fx\leftrightarrow Gx).
    \end{align}
    Then by \cref{metathm:barcan} we have
    \begin{align}
        \forall x^\sigma \ldot \nec(Fx\leftrightarrow Gx).
    \end{align}
    Hence by \labelcref{eq:propositionintensionalism} we have
    \begin{align}
        \forall x^\sigma \ldot (Fx=Gx),
    \end{align}
    which implies \(F=G\) by \labelcref{rule:functionextensionality}.
\end{proof}

\begin{metatheorem}[The Refutation of the Axioms of Extensionality]\label{metathm:nonextensionality}
    The denial of the axiom of propositional extensionality, and hence also of each of the axioms of extensionality for properties and relations, is provable from rules \labelcref{rule:structural,rule:beta,rule:universalinstantiation,rule:universalgeneralization,rule:negationelimination,rule:potentialinfinity}.
    \begin{align}
        \exists pq\ldot ((p\leftrightarrow q)\wedge p\neq q)
    \end{align}
\end{metatheorem}
\begin{proof}
    Suppose otherwise, i.e., \(\forall pq\ldot ((p\leftrightarrow q)\rightarrow p=q)\). Then by \labelcref{rule:structural,rule:beta,rule:universalinstantiation,rule:universalgeneralization,rule:negationelimination} we have
    \begin{align}
        \forall p\ldot  ((p\rightarrow p=\top)\wedge (\neg p\rightarrow p=\bot)),\label{eq:2props}
    \end{align}
    and hence
    \begin{align}
        \neg \exists (1+1+1_t)
    \end{align}
    by Defs.\ \labelcref{eq:def0,eq:def1,eq:def+}. Appealing to \labelcref{eq:2props} again, we have 
    \begin{align}
        \bot = \exists (1+1+1_t)
    \end{align}
    which is contradictory by \labelcref{rule:potentialinfinity}.
\end{proof}

\begin{metatheorem}[Class Comprehension and Class Extensionality; \cite{Church1940-CHUAFO}]
    Where a class is a \(\{\top,\bot\}\)-valued function, i.e.,
    \begin{align}
        \class_\sigma \rewrite \lambda X^{\sigma t}\ldot \forall y\ldot (Xy=\top \vee Xy=\bot),
    \end{align}
    \(\lf\) proves (a) that every property is coextensive with some class:
    \begin{align}
        \forall X^{\sigma t}\ldot \exists Y\in{\class_\sigma}\ldot X\equiv Y,
    \end{align}
    and (b) that classes are extensional:
    \begin{align}
        \forall XY\in{\class_\sigma}\ldot (X\equiv Y\rightarrow X=Y).
    \end{align}
    
\end{metatheorem}
\begin{proof}
    Class comprehension is immediate from Choice; in \(\lf_\iota\) the class coextensive with \(X\) is
    \begin{align}
        \lambda y^{\sigma}\ldot \iota \lambda p\ldot ((Xy\rightarrow p=\top)\wedge (\neg Xy\rightarrow p=\bot)).
    \end{align}
    \parencite[Cf.][p.\ 61.]{Church1940-CHUAFO} The extensionality of classes is immediate from Function Extensionality.
\end{proof}

\begin{metatheorem}[The Necessity of Logic]\label{metathm:noncontingency}
    Where \(\mv{P}\) is a sentence (i.e., a closed formula),
    \begin{align}
        \nec \mv{P}\vee\nec\neg \mv{P}\label{Rule-Necessitation}
    \end{align}
    is provable from rules \labelcref{rule:structural,rule:beta,rule:universalinstantiation,rule:universalgeneralization,rule:negationelimination,rule:Intensionality,rule:functionextensionality,rule:Choice}. Note this no longer holds when \(\mv{P}\) is permitted to contain non-logical constants (including \(\iota\) and \(\varepsilon\)).
\end{metatheorem}
\begin{proof}
    Omitted.
\end{proof}

\begin{metatheorem}[The Complete Atomic Boolean Algebra of Propositions]\label{metathm:completeatomic}
    It is a theorem of \(\lf\) that the propositions form a complete atomic Boolean algebra under \(\wedge\), \(\vee\), \(\neg\). (Hence, by \cref{metathm:actualinfinity}, \(\lf\) proves there are at least \(2^{\aleph_0}\) propositions.)
\end{metatheorem}
\begin{proof}
    Omitted. See Theorem 11.5 of \textcite{Gallin1975-GALIAH}.
\end{proof}

\begin{metatheorem}[The Axioms of Infinity]\label{metathm:actualinfinity}
    Every instance of 
    \begin{align}
        \forall n\in{\nat_\sigma}\ldot \exists n,\label{eq:actualinfinity}
    \end{align}
    is provable in \(\lf\). \labelcref{eq:actualinfinity} asserts that there is an actual infinity of things. 
\end{metatheorem}
\begin{proof}
    Using \cref{metathm:barcan}.
\end{proof}

\begin{metatheorem}[Peano Arithmetic]
    Every theorem of Peano arithmetic is a theorem of \(\lf-\mathrm{\labelcref{rule:Choice}}\) when the constants of arithmetic are defined as above (at any type), and where multiplication is defined in a standard way (omitted here).
\end{metatheorem}
\begin{proof}
    Omitted. 
\end{proof}

\begin{remark}
    \textcite{Goodsell2022-GOOAID-3} shows adding the necessitation of second-order Peano arithmetic for \(\nat_e\) and \(\nat_t\) to the rules \labelcref{rule:structural,rule:beta,rule:universalinstantiation,rule:universalgeneralization,rule:negationelimination,rule:Intensionality,rule:functionextensionality} is equivalent to adding \labelcref{rule:potentialinfinity}.
\end{remark}

\begin{metatheorem}[The Closure of Necessitated Extensions Under the Rule of Intensionality]
    For any set of formulae \(X\), if every formula in \(X\) begins with a \(\nec\) taking widest possible scope, then \(\lf+X\) is closed under the rule of Intensionality (\labelcref{rule:Intensionality}).
\end{metatheorem}
\begin{proof}
    Straightforward from \cref{metathm:intensionalism}.
\end{proof}

\begin{metatheorem}[\(\lf_{\iota}\) and \(\lf_{\varepsilon}\) as Conservative Extensions]
    \(\lf\) is conservatively extended by \(\lf_{\iota}\) and \(\lf_{\varepsilon}\).
\end{metatheorem}
\begin{proof}
    Choice (\labelcref{rule:Choice}) can be used to show, for any finite set of axioms from \(\mathrm{D}_\iota\) or \(\mathrm{C}_\varepsilon\), that functions \(\iota\) or \(\varepsilon\) satisfying those axioms exist. 
\end{proof}

\begin{metatheorem}[Consistency Relative to $\zfc$; Equiconsistency with Henkin's System]
    \(\lf\) is equiconsistent with the system of \textcite{Henkin1950-HENCIT}, and hence is consistent relative to \(\zfc\) by Henkin's soundness theorem.
\end{metatheorem}
\begin{proof}
    Using Benzm{\"u}ller's \parencite*{benzmueller2010simple} method of embedding.
\end{proof}

\section{Relation to other systems}

We now sketch some relations of \(\lf\) to other systems formulated in simply typed languages with abstraction and application, notably omitting various importantly related systems formulated in different languages, such as the aforementioned maximal system $\mathsf{GM}$ that can be constructed out of the components studied in Gallin's \parencite*{Gallin1975-GALIAH} Montague-inspired monograph and Church's \parencite*{church1951formulation} Logic of Sense and Denotation, Alternative (2).

The system of \textcite{Church1940-CHUAFO} is equivalent to \(\lf_\varepsilon\) minus Intensionality (\labelcref{rule:Intensionality}), and with Potential Infinity (\labelcref{rule:potentialinfinity}) replaced with a rule of actual infinity for \(\nat_e\) (but not \(\nat_t\)), i.e., with:
\begin{align}
    \begin{prooftree}
        \hypo{\Gamma\proves \nat_e\mv{n}}
        \infer1{\Gamma\proves \exists_{et}\mv{n}}
    \end{prooftree}\label{eqrule:actualinfinity}
\end{align}

The system of \textcite{Henkin1950-HENCIT} is equivalent to Church's system plus a rule of \textit{extensionality} (notice that, in contrast with \labelcref{rule:Intensionality}, below a \(\Gamma\) occurs on the left-hand side of the turnstile):
\begin{align}
    \begin{prooftree}
        \hypo{\Gamma,\mv{P}\proves \mv{Q}}
        \hypo{\Gamma,\mv{Q}\proves\mv{P}}
        \infer2{\Gamma\proves \mv{P}=\mv{Q}}
    \end{prooftree}\label{eqrule:extensionality}
\end{align}

The system \(\mathsf{HFE}\) of \textcite{Bacon2018-BACTBN} is equivalent to \(\lf\) minus Choice and Potential Infinity (rules \labelcref{rule:Choice,rule:potentialinfinity}). In subsequent work, \textcite{BaconForthcoming-BACC-8} use a weaker system, \textit{Classicism}, which also omits Function Extensionality (\labelcref{rule:functionextensionality}), and strengthens the rule of Intensionality in order to obtain a system that is intensional in the sense of \cref{metathm:fullIntensionality}. A sufficient such strengthening is the following rule, where \(\mv{R}\) and \(\mv{S}\) differ by the substitution of an occurrence of \(\mv{P}\) for \(\mv{Q}\) in any context (including substitutions with variable capture):
\begin{align}
    \begin{prooftree}
        \hypo{\mv{P}\proves \mv{Q}}
        \hypo{\mv{Q}\proves \mv{P}}
        \infer2{\mv{R}\proves \mv{S}}
    \end{prooftree}\label{eqrule:fullIntensionality}
\end{align}

Another axiomatization Bacon and Dorr consider, which is equivalent in a language restricted to relational types (i.e., types that, when angle brackets are omitted in accordance with the convention, end with ``\(t\)'') leaves the rule of Intensionality as it is, but replaces Function Extensionality with the rule
\begin{align}
    \begin{prooftree}
        \hypo{\Gamma\proves \nec\forall \mv{x}\ldot \mv{f}\mv{x}=\mv{g}\mv{x}}
        \infer1{\Gamma\proves \mv{f}=\mv{g}}
    \end{prooftree}\label{eqrule:modalizedfunctionality}
\end{align}
where, as in Function Extensionality, \(\mv{x}\) must not be among the free variables of \(\mv{f}\) or \(\mv{g}\). They show, moreover, that replacing Function Extensionality with this rule in the background of \(\lf\)'s other rules is equivalent with \(\lf\).

\printbibliography

\end{document}